\documentclass[a4paper,11pt,leqno]{article}

\usepackage[latin1]{inputenc}
\usepackage[T1]{fontenc} 
\usepackage[english]{babel}
\usepackage{verbatim}
\usepackage{calligra}
\usepackage{enumerate}
\usepackage{dsfont}
\usepackage{amsfonts}
\usepackage{amsmath}
\usepackage{amsthm}
\usepackage{amssymb}
\usepackage{mathtools}
\usepackage{mathrsfs}
\usepackage{color}
\usepackage{graphicx}
\usepackage{bm}

\usepackage{tikz}
\usepackage[retainorgcmds]{IEEEtrantools}

\usepackage{graphicx}		
\usepackage[hypcap=false]{caption}

\theoremstyle{definition}
\newtheorem{defi}{Definition}[section]

\theoremstyle{plane}
\newtheorem{thm}[defi]{Theorem}

\newcommand{\tsl}{\textsl}

\newcommand{\mc}{\mathcal}

\newcommand{\veps}{\varepsilon}

\newcommand{\what}{\widehat}
\newcommand{\wtilde}{\widetilde}

\newcommand{\ra}{\rightarrow}

\newcommand{\g}{\gamma}

\newcommand{\s}{\sigma}
\renewcommand{\t}{\tau}

\newcommand{\de}{\delta}

\newcommand{\R}{\mathbb{R}}

\newcommand{\N}{\mathbb{N}}

\renewcommand{\Re}{{\rm Re}\,}

\renewcommand{\div}{{\rm div}\,}

\newcommand{\Id}{{\rm Id}\,}

\allowdisplaybreaks

\def\d{\partial}
\def\div{{\rm div}\,}

\begin{document}

\newcommand{\fra}[1]{\textcolor{blue}{FRA: #1}}

\title{\textsc{\Large{\textbf{Well-posedness results for hyperbolic operators \\
with coefficients rapidly oscillating in time}}}}

\author{\normalsize\textsl{Ferruccio Colombini}$\,^1\,,\qquad$ \textsl{Daniele Del Santo}$\,^2\qquad$ and $\qquad$
\textsl{Francesco Fanelli}$\,^{3}$ \vspace{.3cm} \\
\footnotesize{$\,^{1}\;$ \textsc{Universit\`a di Pisa}}  \\
{\footnotesize \it Dipartimento di Matematica} \\
{\footnotesize Largo Bruno Pontecorvo, 5 -- I-56127 Pisa, ITALY} \vspace{.2cm} \\
\footnotesize{$\,^{2}\;$ \textsc{Universit\`a degli Studi di Trieste}} \\
{\footnotesize \it Dipartimento di Matematica e Geoscienze} \\
{\footnotesize Via Valerio, 12/1 -- I-34127 Trieste, ITALY} \vspace{.2cm} \\
\footnotesize{$\,^{3}\;$ \textsc{Universit\'e de Lyon, Universit\'e Claude Bernard Lyon 1}} \\
{\footnotesize \it Institut Camille Jordan -- UMR 5208} \\
{\footnotesize 43 blvd. du 11 novembre 1918, F-69622 Villeurbanne cedex, FRANCE} \vspace{.3cm} \\
\footnotesize{Email addresses: $\,^{1}\;$\ttfamily{ferruccio.colombini@unipi.it}}, $\;$ \footnotesize{$\,^{2}\;$\ttfamily{delsanto@units.it}}, $\;$
\footnotesize{$\,^{3}\;$\ttfamily{fanelli@math.univ-lyon1.fr}}
\vspace{.2cm}
}

\date\today

\maketitle

\subsubsection*{Abstract}
{\footnotesize In the present paper, we consider second order strictly hyperbolic linear operators of the form
$Lu\,=\,\d_t^2u\,-\,\div\big(A(t,x)\nabla u\big)$,
for $(t,x)\in[0,T]\times\R^n$.
We assume the coefficients of the matrix $A(t,x)$ to be smooth in time on $\,]0,T]\times\R^n$, but rapidly oscillating when $t\ra0^+$;
they match instead minimal regularity assumptions (either Lipschitz or log-Lipschitz regularity conditions) with respect to the space variable.

Correspondingly, we prove well-posedness results for the Cauchy problem related to $L$, either with no loss of derivatives (in the Lipschitz case)
or with a finite loss of derivatives, which is linearly increasing in time (in the log-Lipschitz case).

}

\paragraph*{\small 2020 Mathematics Subject Classification:}{\footnotesize 35L10
(primary);
35B65, 35B30, 35L15
(secondary).}

\paragraph*{\small Keywords: }{\footnotesize second order hyperbolic equations; variable coefficients; fast time oscillations; Lipschitz and log-Lipschitz regularity.}


\section{Introduction} \label{s:intro}

The goal of the present paper is to study the well-posedness of the Cauchy problem related to the second order strictly hyperbolic operator
\begin{equation} \label{hyp.op}
Lu\,:=\,\d_t^2u\,-\,\div\big(A(t,x)\,\nabla u\big)\,,\qquad\qquad\qquad \mbox{ with }\qquad \lambda_0\Id\leq A(t,x) \leq \Lambda_0 \Id\,,
\end{equation}
where $\lambda_0$ and $\Lambda_0$ are two positive constants and $A(t,x)\,=\,\big(a_{jk}(t,x)\big)_{j,k=1}^n$ is
defined on a strip $[0,T]\times\R^n$, for some time $T>0$, and is symmetric, \tsl{i.e.} $a_{jk}=a_{kj}$ for any $j,k=1\ldots n$.

It has been well-known since the pioneering work \cite{C-DG-S} (see also \cite {CJS1} for a similar result in the case of weakly hyperbolic equations) that the low regularity in time of the coefficient matrix $A(t,x)$ entails weak well-posedness results,
in the sense that a (finite or even infinte) loss of regularity of the solution $u$ is produced in the time evolution as soon as the coefficients are not Lipschitz-continuous with respect to time. There is a broad literature on the subject,
keeping into account also minimal regularity assumptions on the coefficients with respect to the space variable $x$.
We refer, for instance, to \cite{Nishi}, \cite{Col-Ler}, \cite{Cico-Col}, \cite{DS-Pri1}, \cite{Col-Nishi}, \cite{Col-Kaj}, \cite{T}, \cite{Col-Met} and, more recently, to \cite{CDSFM1}, \cite{CDSFM2}, \cite{Gar-Ruz}, \cite{Cico-Lor}, \cite{CDSF}.
Similar results for first order hyperbolic systems have also been obtained, see \tsl{e.g.} \cite{CDSFM3}, \cite{CDSFM4}.

On the other hand, when dealing with the construction of counterexamples, it appeared from the beginning that phenomena like loss of derivatives, non-existence
or even non-uniqueness of solutions in hyperbolic Cauchy problems may arise when the coefficients are smooth with respect to the time variable, but they
highly oscillate when $t$ approaches one point, which we can assume to be $0$ without loss of any generality. We refer to
\tsl{e.g.} \cite{C-DG-S}, 
\cite{CJS1}, \cite{Col-Spa2}, \cite{CJS2} for the construction of such counterexamples.

As a consequence, some works have focused on studying the relation between fast oscillations of the coefficients $a_{jk}$ and the well-posedness of the Cauchy problem.
In this direction, Yamazaki proved in \cite{Y} that the Cauchy problem for $L$, defined in \eqref{hyp.op}, is well-posed without any loss of derivatives under the following assumptions on the coefficients $a_{jk}$:
they must be continuous in time on $[0,T]\times\R^n$ and $C^2$ in time
on $\,]0,T]\times\R^n$, smooth with respect to the space variable together with $\d_ta_{jk}$ and $\d_t^2a_{jk}$, and they have to satisfy the bound
\begin{equation} \label{hyp:Yam}
\left|t\,\d_ta_{jk}(t,x)\right|\,+\,\left|t^2\,\d_t^2a_{jk}(t,x)\right|\,\leq\,C\qquad\qquad \mbox{ on }\qquad \,]0,T]\times\R^n\,.
\end{equation}
The control on the oscillations of the second order time-derivative seems to be necessary for well-posedness without loss of regularity. As a matter of fact,
for coefficients $a_{jk}=a_{jk}(t)$ depending only on time, in \cite{CDSK1} it was shown that, under the assumption
\[
\left|t\, \frac{d}{dt}a_{jk}(t)\right|\,\leq\,C\qquad\qquad \mbox{ on }\qquad \,]0,T]\times\R^n\,,
\]
the Cauchy problem for $L$ is well-posed with a finite loss of derivatives. In addition, an analogous result was proven in \cite{CDSR}
under the weaker condition
\begin{equation} \label{hyp:faster}
\left|\frac {t}{\log t}\, \frac{d}{dt}a_{jk}(t)\right|\, +\,
\left|\left(\frac {t}{\log t}\right)^2\, \frac{d^2}{dt^2}a_{jk}(t)\right|\,\leq\,C\qquad\qquad \mbox{ on }\qquad \,]0,T]\times\R^n\,.
\end{equation}
We refer also to \cite{R}, \cite{DSKR} and \cite{DS-Pri2} for related results.

\medbreak
In the present paper, we continue the study in this direction and we consider the case of coefficients rapidly oscillating in time, for $t\to0^+$. More preciseley, we will work
under the assumption \eqref{hyp:Yam}, which is the same one formulated in \cite{Y}. Our main contribution
to this problem is two-fold. First of all, we recover Yamazaki's result \cite{Y} by proposing a different proof, which
looks somehow simpler for us (we will give more details about it in a while) than the original one.
Secondly, and more importantly, we will be able to consider minimal regularity assumptions
on the coefficients with respect to the space variable, namely Lipschitz and log-Lipschitz regularities in $x$. In the former case (Lipschitz regularity),
we will prove an energy estimate for $L$ with \emph{no loss} of derivatives in the space $H^{1-\theta}\times H^{-\theta}$, for any $\theta\in[0,1[\,$. From this, it
is classical to recover well-posedness of the Cauchy problem related to $L$ in the same space. In the latter case (log-Lipschitz regularity of
the coefficients in $x$), we will prove an energy estimate with a finite, linearly increasing in time, \emph{loss of derivatives}:
if the datum $\big(u_0,u_1\big)$ belongs to $H^{1-\theta}\times H^{-\theta}$, for some $0<\theta<1$, then the norm of the solution at time $t$
is controlled in the space $H^{1-\theta-\beta^*t}\times H^{-\theta-\beta^*t}$, where the constant $\beta^*>0$ depends only on the coefficient matrix $A(t,x)$.

The proof is inspired by the techniques recently employed in \cite{CDSFM1}, \cite{CDSFM2} and \cite{CDSF}, and consists of two main ingredients. The first ingredient
is a definition of the energy functional which presents a lower order corrector; in turn, this idea goes back to work \cite{T} of Tarama. The role of the
corrector is to cancel out bad remainders which appear when computing the time derivative of the energy and which would be otherwise out of control.
The second ingredient, instead, is the use of paradifferential calculus with parameters, as introduced by M\'etivier (see \cite{Met_1986} and \cite{Met-Z}),
which provides us with a powerful tool for dealing with the low regularity in space of the coefficients.

To conclude this short introduction, let us mention that, in the case when $a_{jk}=a_{jk}(t)$ does not depend on $x$,
we are able to consider even faster oscillations, in the spirit of condition \eqref{hyp:faster} from \cite{CDSR},
recovering the result of \cite{CDSR} with a more elementary proof. For this, it is enough to use exactly the same energy of Tarama \cite{T} in order to carry out
energy estimates.
However, at present we are not able to extend this approach to the case in which the coefficients also depend on the space variables. This is why we will not devote attention
to that case in this work.

\subsubsection*{Acknowledgements}

{\small
The first and the second author are members of the Gruppo Nazionale per l'Analisi Matematica, la Probabilità e le
loro Applicazioni (GNAMPA) of the Istituto Nazionale di Alta Matematica (INdAM). They acknowledge INdAM for the support.

The work of the third author has been partially supported by the LABEX MILYON (ANR-10-LABX-0070) of Universit\'e de Lyon, within the program
``Investissement d'Avenir'' (ANR-11-IDEX-0007),  and by the projects SingFlows (ANR-18-CE40-0027)
and CRISIS (ANR-20-CE40-0020-01), both operated by the French National Research Agency (ANR).

}

\section{Statement of the main results} \label{s:results}

We consider the linear differential operator
\begin{equation}\label{def:L}
Lu\,=\, \partial_t^2 u\,-\,\sum_{j,k=1}^n \partial_j\Big(a_{jk}(t,x)\partial_k u\Big)\,,
\end{equation}
where $a_{j,k}:[0,T] \times \R^n\to \R$ are such that $a_{jk}=a_{kj}$.
We also assume that  there exist $\lambda_0$, $\Lambda_0>0$ such that, for all $(t,x)\in [0,T] \times \R^n$ and $\xi\in \R^n$, 
\begin{equation} \label{hyp:hyp}
\lambda_0 |\xi|^2\leq \sum_{j,k=1}^n a_{jk}(t,x)\xi_j\xi_k \leq \Lambda_0|\xi|^2. 
\end{equation}
In particular, under this condition, the operator $L$ becomes strictly hyperbolic.

Suppose moreover that the coefficients $a_{j,k}$ satisfy the following assumptions.
\begin{itemize}
\item[\bf (H1)] 
There exist $\ell\in\{0,1\}$ and $C_0>0$ such that, for any $j,k\in\{1\ldots n\}$, one has
$$
|a_{jk}(t, x+y)-a_{jk}(t, x)|\leq C_0 |y| \log^\ell\left(1+ \frac{1}{|y|}\right)\,,
$$
for all $t\in [0,T]$ and for all $x,y\in \R^n$, with $0<|y|<1$.

\item[\bf (H2)] Each $a_{jk}$ is twice differentiable with respect to $t$ in  $\,]0,T] \times \R^n$; in addition
there exist constants $C_1>0$, $C_2>0$, $C_3>0$ for which, for any $j,k\in\{1\ldots n\}$, the controls
\begin{align*}
\left|\partial_t a_{jk}(t,x)\right|&\leq C_1\,  \frac{1}{t}\,,\\
\left|\partial_t  a_{jk}(t,x+y)- \partial_t a_{jk}(t,x)\right|&\leq C_2\,  \frac{1}{t}\,  |y| \log^\ell \left(1+ \frac{1}{|y|}\right)\,,\\
\left|\partial^2_t a_{jk}(t,x)\right|&\leq C_3 \, \frac{1}{t^2}
\end{align*}
hold true for all $t\in \, ]0,T]$ and for all $x,y\in \R^n$, with $0<|y|<1$, where $\ell\in\{0,1\}$ is the same exponent fixed in assumption (H1).
\end{itemize}

Notice that hypothesis (H1) means that the coefficients $a_{jk}$ are Lipschitz continuous in space when $\ell=0$, log-Lipschitz
continuous when $\ell=1$, uniformly with respect to the time variable. Requiring a similar control also for $\d_ta_{jk}$, see the second condition appearing in (H2),
is somehow natural in this kind of study; see \tsl{e.g.} \cite{CDSFM1}, \cite{CDSFM2}, \cite{CDSF}.

\medbreak
We are now able to state the main results of the paper.
The first one is devoted to the Lipschitz case, namely to the case $\ell=0$. It is an extension of Yamazaki's result \cite{Y} to the case when
coefficients have low regularity in the space variable; besides, we propose a proof based on a different technique, which looks to an extent simpler.

\begin{thm} \label{th:Yam}
Let $L$ be defined by \eqref{def:L}, and assume that the assumptions \eqref{hyp:hyp}, (H1) and (H2) are in force with $\ell=0$.
Fix $\theta\in[0,1[\,$.

Then, there exists a universal constant $C>0$ such that, for any $u\in C^2\big([0,T]; H^\infty(\R^n)\big)$, one has the estimate
\begin{equation}\label{noloss}
\begin{array}{ll}
\displaystyle{\sup_{0\leq t \leq T}\Big(\|u(t,\cdot)\|_{H^{1-\theta}}+ \|\partial_t u(t,\cdot)\|_{H^{-\theta}}\Big)}\\[0.2 cm]
\qquad\qquad\qquad\displaystyle{
\leq C\,\left(\|u(0,\cdot)\|_{H^{1-\theta}}+ \|\partial_t u(0,\cdot)\|_{H^{-\theta}}+\int_0^{T}\|Lu(\t, \cdot)\|_{H^{-\theta}}\, d\t\right)\,.}
\end{array}
\end{equation}
\end{thm}

From the previous statement, it is classical to derive the well-posedness (without loss of derivatives) of the Cauchy problem associated to $L$
in the energy space $H^{1-\theta}(\R^n)\times H^{-\theta}(\R^n)$, for any given $\theta\in[0,1[\,$.

\medbreak
Let us now switch to the log-Lipschitz case, for which $\ell=1$ in (H1) and (H2) above. As is well-known (se \tsl{e.g.}
\cite{Col-Ler}, \cite{Col-Met}, \cite{CDSFM1} and references therein), in this case the estimates present a loss of derivatives, which is linearly increasing in time.
In addition, because of that and because of the low regularity in space of the coefficients, the estimate becomes local in time. Finally, one has to work in energy
spaces which are less regular than $H^1(\R^n)\times L^2(\R^n)$, \tsl{i.e.} one has to take $\theta\in\,]0,1[\,$.

\begin{thm}\label{theorem 3}
Let $L$ be defined by \eqref{def:L}, and assume that its coefficients $a_{jk}$ satisfy assumptions \eqref{hyp:hyp}, (H1) and (H2) with $\ell=1$.
Fix $\theta\in \, ]0, 1[\,$.

Then, there exist $\beta^*>0$, $T_0 \in \,] 0,T]$ and  $C>0$ such that, for all function $u\in C^2\big([0,T_0]; H^\infty(\R^n)\big)$, one has
the energy estimate
\begin{equation}\label{lossdept}
\begin{array}{ll}
\displaystyle{\sup_{0\leq t \leq T_0}\Big(\|u(t,\cdot)\|_{H^{1-\theta-\beta^* t}}+ \|\partial_t u(t,\cdot)\|_{H^{-\theta -\beta^* t}}\Big)}\\[0.2 cm]
\qquad\qquad\displaystyle{
\leq C\left(\|u(0,\cdot)\|_{H^{1-\theta}}+ \|\partial_t u(0,\cdot)\|_{H^{-\theta}}+\int_0^{T_0}\|Lu(\t, \cdot)\|_{H^{-\theta-\beta^*\t}(\R^n)}\, d\t\right)\,.}
\end{array}
\end{equation}
\end{thm}

From the previous result, the well-posedness of the Cauchy problem related to $L$ follows in the space $H^\infty$, with a finite loss of derivatives.

\medbreak
The rest of the paper is devoted to the proof of the previous results. As already announced, we propose a different approach than \cite{Y}, based on
Tarama's energy \cite{T} and paradifferential calculus with parameters, as used in \cite{CDSFM1} and \cite{CDSFM2} for operators with low regularity coefficients
both in time and space variables.

Before moving on, we remark that, in the case that the coefficients $a_{jk}=a_{jk}(t)$ depend only on time, our method enable us to treat the case
of even faster oscillations in time. Namely, in assumption (H2) we can admit an additional logarithmic loss: there exists $\de\in[0,1]$ such that
\begin{align*}
\left|\partial_t a_{jk}(t)\right|&\leq C\,  \frac{1+\left|\log t\right|^\delta}{t}\,,\\
\left|\partial^2_t a_{jk}(t)\right|&\leq C \,\left(\frac{1+\left|\log t\right|^\delta}{t}\right)^2\,.
\end{align*}
By resorting to Tarama's energy (defined in frequencies, as in \cite{C-DG-S}, \cite{T}), we can prove energy estimates with no loss when $\de=0$, with finite
(fixed) loss of derivatives when $\de>0$ and with arbitrarily small loss when $0<\de<1$, in any Sobolev space $H^s$. This result has already been known
(see \cite{CDSR}, \cite {R} for the case $\de=1$ and \cite{Y} for the case $\de=0$), but our proof looks much simpler than the previous ones.

\section{Approximation and regularisation} \label{s:approx}

The goal of this section is to introduce a suitable approximation/regularisation of the coefficients, both in time and space variables, in order
to perform (in Section \ref{s:energy}) energy estimates for approximated energies related to the function $u$.

\subsection{Approximation of the coefficients} \label{ss:time}

We start by dealing with the time variable. In order to perform energy estimates, we need to define approximations of the coefficients which
are smooth up to the time $t=0$.
The obvious problem is that, under assumption (H2), the coefficients are highly oscillating  when $t$ approaches $0$. As a consequence, we need to approximate
the original coefficients by smooth ones: for this, we will follow \cite{CDSR} and proceed by truncation around $t=0$.

Thus, let  $\varepsilon\in \,\left] 0,\frac {T}{2}\right]$. Following \cite{CDSR}, we set
$$
{\wtilde a}_{\varepsilon, jk}(t,x)\,:=\,\left\{
\begin{array}{ll}
a_{jk}(\varepsilon, x)\qquad&\qquad\text{for} \quad t<\varepsilon\,,\\[0.2 cm]
a_{jk}(t, x)\qquad&\qquad\text{for} \quad \varepsilon\leq t\leq T\,,\\[0.2 cm]
a_{jk}(T, x)\qquad&\qquad\text{for} \quad t>T\,.
\end{array}
\right.
$$
Of course, $\wtilde a_{jk,\veps}$ previously defined is not smooth in time, so not suitable for energy estimates. Therefore,
a more involved approximation procedure is required.
For this, consider $\big(\rho_\eta\big)_{\eta\in\,]0,1]}$ a standard family of mollifiers. We then define 
$$
{\wtilde{\wtilde a}}_{\varepsilon, jk}(t,x)\,:=\, \big(\rho_{\varepsilon/2}\,*_t\, {\wtilde a}_{\varepsilon, jk}\big)(t,x)\,,
$$
where the notation $*_t$ means that the convolution is taken only with respect to the time variable.

This approximation is not satisfactory for our scopes yet, because, in view of the next computations, it would be more convenient to
keep the old values of the function $a_{jk}$ when the time $t$ is far from $0$.
Hence, let us take a function
$\theta \in C^\infty(\R)$ such that $0\leq \theta\leq 1$, $\theta(s)=1$ for $s\leq 1/4$ and  $\theta (s)=0$ for $s\geq 3/4$. Then,
we set $\theta_\varepsilon(t)\,:=\,\theta\left(\frac{1}{3\varepsilon}t+\frac{1}{12}\right)$ and we define
\begin{equation}\label{approx}
a_{\varepsilon, jk}(t,x)\,:=\,{\wtilde{\wtilde a}}_{\varepsilon, jk}(t,x)\,\theta_\varepsilon(t)\,+\, a_{jk}(t, x)\,\big(1-\theta_\varepsilon(t)\big)\,.
\end{equation}
We remark that, with this definition, one has
$$
\begin{array}{ll}
\displaystyle{a_{\varepsilon, jk}(t,x)= a_{jk}(t,x)}\qquad &\text{ if }\quad 2\veps\leq t\leq T\,,\\[0.2 cm]
\displaystyle{a_{\varepsilon, jk}(t,x)= a_{jk}(\varepsilon,x)}\qquad &\text{ if }\quad \displaystyle{t\leq \frac{\varepsilon}{2}}
\end{array}
$$
for any $x\in\R^n$. In addition, the coefficients $a_{\veps,jk}$ satisfy the following properties:
\begin{enumerate}[(i)]
\item for all $(t,x)\in \R \times \R^n$ and $\xi\in \R^n$, we have
$$
\lambda_0 |\xi|^2\leq \sum_{j,k=1}^n a_{\varepsilon, jk}(t,x)\xi_j\xi_k \leq \Lambda_0|\xi|^2\,; 
$$
\item there exists a constant $C_0>0$, possibly different from the one appearing in (H1) in Section \ref{s:results}, but still independent of $\veps\in\,]0,1]$, such that 
$$
\left|a_{\varepsilon, jk}(t,x+y)- a_{\varepsilon, jk}(t,x)\right|\leq C_0\, |y|\, \log^\ell\left(1+ \frac{1}{|y|}\right)
$$
for all $t\in \, {\R}$ and for all $x,y\in \R^n$, with $0<|y|<1$;
\item there exists a function $\varphi_\varepsilon : \R\to\R $ of class $C^2$, satisfying
$$
\varphi_\varepsilon(t)=0 \quad\text{ for }\quad t\leq \frac{\varepsilon}{2} \qquad\qquad\text{ and }\qquad\qquad
0\leq \varphi_\varepsilon(t)\leq \frac{1}{t} \quad\text{ for }\quad t\geq \frac{\varepsilon}{2}\,,
$$
and there exist constants $C_1$, $C_2$, $C_3>0$ (possibly different from the ones appearing in (H2) of Section \ref{s:results}, but still independent
of $\veps\in\,]0,1]$) such that
\begin{align*}
\left|\partial_t a_{\varepsilon, jk}(t,x)\right|&\leq C_1\,\varphi_\varepsilon(t)\,, \\
\left| \partial_t a_{\varepsilon, jk}(t,x+y)- \partial_t a_{\varepsilon, jk}(t,x) \right|&\leq C_2\,  \varphi_\varepsilon(t)\, |y|\,
\log^\ell\left(1+ \frac{1}{|y|}\right)\,,\\
\left|\partial^2_t a_{\varepsilon, jk}(t,x)\right|&\leq C_3\, \big(\varphi_\varepsilon(t)\big)^2\,,
\end{align*}
for all $t\in \R$ and for all $x,y\in \R^n$, with $0<|y|<1$.
\end{enumerate}

The coefficients $a_{\veps,jk}$ constructed in this way are certainly smooth in time on $[0,T]\times\R^n$, but they are not regular with respect to the space variable.
Therefore, a further approximation is needed: this is the scope of the next subsections.

\subsection{Tools from Littlewood-Paley theory} \label{ss:LP}

Here we recall basic facts linked to Littlewood-Paley decomposition and paradifferential calculus with parameters, which will be very useful tools in our
proof.

\subsubsection*{Littlewood-Paley decomposition and Sobolev spaces}

Let us start by defining the Littlewood-Paley decomposition of any tempered distribution $u\in\mc S'(\R^n)$ and the dyadic characterisation of Sobolev spaces.
We refer to \tsl{e.g.} \cite{BCD} and \cite{Met_2008} for details.

We fix a smooth radial function $\chi:\R^n\to\R$, supported in the ball $B(0,2)$ of center $0$ and radius $2$, with $\chi\equiv1$ in a neighbourhood 
of the unit ball $B(0,1)$ and $\chi$ non-increasing along the radial directions. For $\xi\in\R^n$, we also define $\phi(\xi)\,:=\,\chi(\xi)-\chi(2\xi)$.

For any $j\in\N$ and $\xi\in\R^n$, we also define
\begin{equation} \label{def:LP_loc}
\chi_j(\xi)\,:=\,\chi\big(2^{-j}\xi\big)\qquad\qquad \mbox{ and }\qquad\qquad \phi_j(\xi)\,:=\,\phi\big(2^{-j}\xi\big)\,.
\end{equation}

Then, we can introduce the dyadic blocks $\Delta_j$, for $j\in\N$. First of all, we define $\Delta_0$ as the operator
\[
\Delta_0u\,:=\,\chi(D_x)u\,=\,\mc F^{-1}\big(\chi(\xi)\,\what{u}(\xi)\big)\,,
\]
where $D_x=-i\d_x$, $\what u$ denotes the Fourier transform of a tempered distribution $u$ and the symbol $\mc F^{-1}$ the inverse Fourier transform.
Next, for any $j\geq1$, we set
\[
\Delta_ju\,:=\,\phi_j(D_x)u\,=\,\mc F^{-1}\big(\phi_j(\xi)\,\what u(\xi)\big)\,.
\]
For $j\geq0$, we also set
\[
S_ju\,:=\,\chi_j(D_x)u\,=\,\sum_{k=0}^j\Delta_ku\,.
\]


Thanks to the previous operators, one can define the \emph{Littlewood-Paley decomposition} of tempered distributions:
for any $u\in\mc S'(\R^n)$, the decomposition
\[
u\,=\,\sum_{j=0}^{+\infty}\Delta_ju 
\]
holds true in the sense of $\mc S'(\R^n)$.

One can also characterise Sobolev classes in terms of Littlewood-Paley decomposition. More precisely, 
given some $s\in\R$ and $u\in\mc S'(\R^n)$, one has that
\[
u\in H^s(\R^n)\qquad\quad \Longleftrightarrow\qquad\quad \Big(2^{sj}\,\left\|\Delta_ju\right\|_{L^2}\Big)_{j\geq0}\in\ell^2(\N)\,,
\]
with equivalence of norms: there exists a constant $C>0$, independent of $u$ but possibly depending on $s$, such that
\[
\frac{1}{C}\,\sum_{j=0}^{+\infty}2^{2sj}\,\left\|\Delta_ju\right\|_{L^2}^2\,\leq\,\left\|u\right\|_{H^s}^2\,\leq\,
{C}\,\sum_{j=0}^{+\infty}2^{2sj}\,\left\|\Delta_ju\right\|_{L^2}^2\,.
\]

\subsubsection*{Paradifferential calculus with parameters}

In order to deal with the low regularity of the coefficients with respect to the space variable, we will adopt the approach used in \cite{CDSFM1} and \cite{CDSFM2},
based on a paralinearisation of the elliptic part $\sum_{j,k}\d_j\big( a_{jk}\,\d_k\big)$ of the operator $L$.
Of course, $L$ being linear, we only mean that the products by the coefficients $a_{jk}$ will be replaced by paraproduct operators $T_{a_{jk}}$
(see \cite{Met_2008}, \cite{BCD}): as a matter of fact, this procedure intrinsically involves a regularisation in space.

On the other hand, in order to guarantee positivity of the paralinearised operator, we have to resort to M\'etivier's paradifferential calculus with parameters,
as introduced in \cite{Met_1986}, \cite{Met-Z}.
This technique being by now quite well understood, we will just recall the main ideas of the construction and use it in our computations;
we refer the reader to \cite{CDSFM1} and \cite{CDSFM2} for details.

\medbreak
Thus, let us fix a symbol $f=f(x,\xi)$, defined for $(x,\xi)\in\R^n\times\R^n$, smooth with respect to the $\xi$ variable and merely bounded in $x$. The main idea of
paradifferential calculus (see \cite{Met_2008}) is to associate to $f$ a classical symbol $\s_f=\s_f(x,\xi)$ which is smooth also with respect to the space variable
$x$. Then, one defines the paradifferential operator $T_f$ associated to the (rough) symbol $f$ as the classical pseudodifferential operator associated to $\s_f$: we
set
\[
 T_f\,:=\,\s_f(x,D_x)\,.
\]
If the symbol $f$ also depend on time, we can treat the time variable as a parameter in the construction.

This having beeing defined, one can thus develop a symbolic calculus for pseudodifferential operators, based on the smoothness of the original symbol $f$.
We will see in a while how to do that in the case of our hyperbolic operator $L$. For the time being, let us emphasise a couple of facts.

First of all, all the construction is, in some sense, \emph{canonical}: although the classical symbol $\s_f$ associated to $f$ is not uniquely determined,
the paradifferential operator $T_f$ is, up to a remainder of lower order (the precise order depends on the space regularity and frequency decay of $f$).
In particular, one can fix the choice
\[
T_f:\;u\,\mapsto\,T_fu\,:=\,S_{0}f\,S_3u\,+\,\sum_{\nu=1}^{+\infty}S_\nu f\,\Delta_{\nu+3}u\,,
\]
which coincides (whenever $f=f(x)$ only depends on the space variable $x$) with the classical Bony paraproduct operator \cite{Bony} (see also
Chapter 2 of \cite{BCD}). The point is that this construction makes sense for more general symbols $f=f(t,x,\xi)$.

Nonetheless, we remark that, with this construction, for positive symbols one loses positivity of the corresponding paradifferential operator,
a property which is very much desirable when dealing with energy functionals in hyperbolic Cauchy problems.
As a matter of fact, if $f=f(x)$ is a bounded scalar function defined over $\R^n$, with $f(x)\geq c>0$
for a.e. $x\in\R^n$, one has $\left\|f\,u\right\|_{L^2}\,\geq\,c\,\|u\|_{L^2}$, but it is not true, in general,
that the property $\left\|T_fu\right\|_{L^2}\,\geq\,\wtilde{c}\,\|u\|_{L^2}$ holds (even with a new constant $\wtilde c>0$).

In order to solve the previous issue, the basic idea of Métivier (see \cite{Met_1986}, \cite{Met-Z}) was to add
enough blocks of frequencies to the definition of the paraproduct to recover, in turn, positivity of the modified operator. Notice that the product $fu$ mixes all the frequencies of $f$ and $u$, whereas the paraproduct $T_fu$ uses only a few of them. More rigorously,
one introduces a parameter $\g\geq1$ in the construction,
obtaining classical symbols $\s_f^\g$ and operators $T_f^\g$ depending on such $\g$. Considering larger values of $\g$ allows to take into account more
and more frequencies in the definition of the operator $T_f^\g$: for $u\in\mc S(\R^n)$, one can consider
\begin{equation} \label{def:paraprod}
T_f^\g u\,=\,S_{\mu-1}f\,S_{\mu+2}u\,+\,\sum_{\nu=\mu}^{+\infty}S_{\nu}f\,\Delta_{\nu+3} u\,,
\end{equation}
where $\mu=\mu(\g)=\left[\log_2\g\right]$ is the integer part of $\log_2\g$. Now, given a symbol $f=f(x,\xi)$ of (say) order $0$ such that $f\geq c>0$,
there exists a $\g_0\,=\,\g_0(f)$, depending on the symbol $f$, such that, for all $\g\geq \g_0$, one has
\[
\forall\,u\in\mc S(\R^n)\,,\qquad\qquad 
\left\|T_f^\g u\right\|_{L^2}\,\geq\,\frac{c}{2}\,\|u\|_{L^2}\,.
\]

In light of the previous considerations, without loss of generality, throughout this paper we can fix the previous choice
of paradifferential operators: given a symbol $f=f(t,x,\xi)$ as described above, we define the operator $T_f^\g$ according to
formula \eqref{def:paraprod}.
Observe also that, more generally, one can even deal with symbols $f=f(t,x,\xi,\g)$ depending also on the parameter $\g$.

\medbreak
To conclude this part, let us recall (see Section 3.2 of \cite{CDSFM1} for details) that we can define also a Littlewood-Paley decomposition
depending on the parameter $\g\geq1$, as well as Sobolev spaces $H^s_\g$ depending on $\g$.
For $\g\geq1$ fixed, the space $H^s_\g$ coincides with the classical $H^s$ space, with equivalence of norms (where the multiplicative constants depend on $\g$, of course).
Thus, the dyadic characterisation that we have seen above still applies.

\subsection{Symbolic calculus} \label{ss:symbolic}

Let us see how to apply the previous ideas to the case of the strictly hyperbolic operator $L$, defined in \eqref{def:L}.
To begin with, we introduce the $0$-th order symbol
$$
\alpha_\varepsilon(t,x,\xi,\gamma)\,:=\,\left(\gamma^2+|\xi|^2\right)^{-\frac{1}{2}}\,\left(\gamma^2+\sum_{j,k=1}^n a_{\varepsilon, jk}(t,x)\xi_j\xi_k\right)^{\frac{1}{2}}\,.
$$

Following the original idea of \cite{C-DG-S} and later \cite{Col-Ler}, we immediately choose
$$
\varepsilon=2^{-\nu}\,.
$$
For simplicity of notation, we write $\alpha_\nu$, $a_{\nu, jk}$ and $\varphi_\nu$ instead of $\alpha_{2^{-\nu}}$, $a_{2^{-\nu},jk}$ and $\varphi_{2^{-\nu}}$, respectively.

In order to develop symbolic calculus, we need suitable estimates for the classical symbols associated to the non-smooth symbol $\alpha_{\nu}$
and its time derivatives.
Repeating the computations of \cite{CDSFM1}, \cite{CDSFM2} and \cite{CDSF}, it is easy to obtain the following bounds:
for any $\alpha, \beta \in\N^n$, with $|\beta|\geq1$, and for any $(t,x,\xi,\g)\in[0,T]\times\R^n\times\R^n\times[1,+\infty[\,$, one has
\begin{align*}
\left|\partial_\xi^\alpha \sigma_{\alpha_\nu}(t,x,\xi,\gamma)\right| &\leq C_\alpha\,\left(\gamma+|\xi|\right)^{-|\alpha|}\,,\\[0.2 cm]
\left|\partial_x^\beta\partial_\xi^\alpha \sigma_{\alpha_\nu}(t,x,\xi,\gamma)\right| &\leq
C_{\beta, \alpha}\,\left(\gamma+|\xi|\right)^{-|\alpha|+|\beta|-1}\,\log^\ell\left(1+\gamma+|\xi|\right)\,,\\[0.2 cm]
\left|\partial_\xi^\alpha \sigma_{\partial_t \alpha_\nu}(t,x,\xi,\gamma)\right| &\leq
C_{\alpha}\,\varphi_\nu (t)\, \left(\gamma+|\xi|\right)^{-|\alpha|}\,,\\[0.2 cm]
\left|\partial_x^\beta\partial_\xi^\alpha \sigma_{\partial_t \alpha_\nu}(t,x,\xi,\gamma)\right| &\leq
C_{\beta, \alpha}\,\varphi_\nu (t)\,\left(\gamma+|\xi|\right)^{-|\alpha|+|\beta|-1}\,\log^\ell \left(1+\gamma+|\xi|\right)\,,\\[0.2 cm]
\left|\partial_x^\beta\partial_\xi^\alpha \sigma_{\partial^2_t \alpha_\nu}(t,x,\xi,\gamma)\right| &\leq
C_{\alpha}\, \big(\varphi_\nu (t)\big)^2\, \left(\gamma+|\xi|\right)^{-|\alpha|+|\beta|}\,, 
\end{align*}
where the involved multiplicative constants do not depend on $\nu$, nor on $(t,x,\xi,\g)$.
We point out that, whenever an estimate as above holds with $\ell=1$, we will say that the paradifferential operator associated to the symbol under consideration
is of order ${\rm Log}$.

Next, we fix $\gamma\geq 1$ so large that, for all $u\in H^\infty$, one has
\begin{equation}\label{positivity2}
\left\|T^\gamma_{\alpha_{\nu}^{-1/2}} u\right\|_{L^2} \geq \frac{\lambda_0}{2}\,\|u\|_{L^2}\qquad\quad\text{ and }\qquad\quad
\left\|T^\gamma_{\alpha_{\nu}^{1/2}(\gamma^2+|\xi|^2)^{1/2}} u\right\|_{L^2}\geq \frac{\lambda_0}{2}\|u\|_{H^1}\,.
\end{equation}
From now on, we forget about the index $\g$ (which is now fixed) in the notation and we simply write 
$T_{f}$ instead of $T^\gamma_{f}$, for any symbol $f$ which will enter into play in our computations.

\section{Energy estimates} \label{s:energy}

We are now in the position of proving the energy estimates claimed in Theorems \ref{th:Yam} and \ref{theorem 3}. For this, as already announced, we will
follow the strategy of \cite{Col-Ler}, \cite{CDSFM1} and \cite{CDSF}.

To begin with, some preliminaries are in order.
Firstly, we approximate the operator $L$ by paralinearisation: given $u\in C^2\big([0,T];H^\infty(\R^n)\big)$,
we have
\begin{equation} \label{eq:paralin}
\partial_t^2 u\,=\, \sum_{j,k=1}^n\partial_j\big(a_{jk}(t,x)\partial_k u\big)+Lu\,=\,
\sum_{j,k=1}^n\partial_j\left(T_{a_{jk}}\partial_k u\right)+\wtilde Lu\,,
\end{equation}
where we have defined
$$
\wtilde Lu\,:=\,Lu+\sum_{j,k=1}^n\partial_j\Big(\left(a_{jk}-T_{a_{jk}}\right)\partial_k u\Big)\,.
$$

Next, we want to localise the previous equation in frequency, getting in this way an equation for each dyadic block $u_\nu:=\Delta_\nu u$ associated to $u$.
Applying the operator $\Delta_\nu$ to equation \eqref{eq:paralin} yields
\begin{equation} \label{eq:u_nu}
\partial_t^2 u_\nu\,=\, \sum_{j,k=1}^n\partial_j\left(T_{a_{jk}}\partial_k u_\nu\right)
+ \sum_{j,k=1}^n\partial_j\Big(\left[\Delta_\nu, T_{a_{jk}}\right]\partial_k u\Big)+ \left(\wtilde L u\right)_\nu\,,
\end{equation}
where we have set $\left(\wtilde L u\right)_\nu=\Delta_\nu\left(\wtilde L u\right)$ and $\left[\Delta_\nu, T_{a_{jk}}\right]$ denotes the commutator between the
localisation operator $\Delta_\nu$ and the paramultiplication operator $T_{a_{jk}}$.

Our first goal is to estimate, by the use of \eqref{eq:u_nu}, the growth of the localised energy, namely of the energy $e_\nu$ associated to each dyadic block
$u_\nu$. The total energy $\mc E$ will be defined only afterwards, as a suitable weighted sum of the $e_\nu$'s.
Observe that, owing to the fast oscillations of the coefficients in time, defining $e_\nu$ in a classical way, namely as
the sum of $\|\d_tu_\nu\|_{L^2}^2\,+\,\|u_\nu\|^2_{H^1}$, would cause the appearing of bad terms
in the energy estimates, which are out of control. The key point here is to resort to the idea of Tarama \cite{T} and introduce a lower order corrector
in the definition of $e_\nu$, whose role is exactly to cancel those bad terms. We refer to \cite{CDSFM1}, \cite{CDSFM2} and \cite{CDSFM3} for further
implementations of this idea in the context of hyperbolic Cauchy problems.

\subsection{The localised energy} 

In order to use Tarama's energy in our context, we start by defining
\begin{align*}
v_\nu(t,x)\,&:=\, T_{\alpha^{-1/2}_{\nu} }\partial _t u_\nu- T_{\partial_t( \alpha^{-1/2}_{\nu})} u_\nu\,, \\[0.2 cm]
w_\nu(t,x)\,&:=\,T_{\alpha^{1/2}_{\nu}(\gamma^2+|\xi|^2)^{1/2 }}\ u_\nu\,.
\end{align*}
Then, we introduce the approximate energy associated to the $\nu$-th component $u_\nu$ as
$$
e_\nu(t)\,:=\,\left\|v_\nu(t, \cdot)\right\|^2_{L^2}+\left\|w_\nu(t, \cdot)\right\|^2_{L^2}+\left\|u_\nu(t, \cdot)\right\|^2_{L^2}\,.
$$

We remark that $e_\nu$ is equivalent to the classical energy associated to $u_\nu$. As a matter of fact, we can bound
\begin{align*}
\left\|\partial_t u_\nu(t,\cdot)\right\|_{L^2} & \leq C  \left\|T_{ \alpha^{-1/2}_{\nu}} \partial_t  u_\nu(t,\cdot)\right\|_{L^2} \\[0.2 cm]
&\leq C\,\left(\left\|v_\nu(t,\cdot)\right\|_{L^2}+  \left\|T_{\partial_t( \alpha_{\nu}^{-1/2})} u_\nu(t,\cdot)\right\|_{L^2}\right) \\[0.2 cm]
&\leq C\,\Big(\left\|v_\nu(t,\cdot)\right\|_{L^2}+ \varphi_\nu (t)  \left\|u_\nu(t,\cdot)\right\|_{L^2}\Big)\\[0.2 cm]
&\leq C\,\Big(\left\|v_\nu(t,\cdot)\right\|_{L^2}+ \varphi_\nu (t)\,2^{-\nu}  \left\|w_\nu(t,\cdot)\right\|_{L^2}\Big)\,.
\end{align*}
At this point, noticing that
\begin{equation}\label{boundphi2}
0\,\leq\,\varphi_\nu (t)\,2^{-\nu}\,\leq\,C\,,
\end{equation}
for a suitable constant $C>0$, we deduce that 
\begin{equation}\label{estpartunu2}
\left\|\partial_t u_\nu (t, \cdot)\right\|_{L^2}\leq C\,\big(e_\nu(t)\big)^{\frac{1}{2}}\,.
\end{equation}
From this inequality, together with \eqref{positivity2}, we easily deduce that 
\begin{equation}\label{estenu2}
\forall\,t\in[0,T]\,,\qquad\qquad
e_\nu(t)\,\sim\, \left\|\partial_t u_\nu(t,\cdot)\right\|^2_{L^2} + \left\|u_\nu(t,\cdot)\right\|^2_{L^2} + \left\|\nabla_x u_\nu(t,\cdot)\right\|^2_{L^2}\,,
\end{equation}
where the (implicit) multiplicative constants do not depend on $t\in\R_+$, nor on $u$ and $\nu\in\N$.

\subsection{Time derivative of the energy components}

We now proceed to differentiating $e_\nu$ with respect to time and estimating the terms thus produced.

\subsubsection{The zero-th order term}

To begin with, we compute the time derivative of the zero-th order term appearing in the definition of $e_\nu$. This is an easy task: we have
\[
\frac{d}{dt}\|u_\nu(t,\cdot)\|_{L^2}^2\,=\, 2\, \Re\big(u_\nu(t,\cdot), \partial_t u_\nu(t,\cdot)\big)_{L^2}\,.
\]
Consequently, in view of the bound \eqref{estpartunu2}, we obtain
\begin{equation}\label{estdtunu2}
\frac{d}{dt}\|u_\nu(t,\cdot)\|_{L^2}^2\, \leq\, C\, e_\nu(t)\,.
\end{equation}

The computations of the terms involving $v_\nu$ and $w_\nu$ are a bit more involved: we now switch our attention to them.

\subsubsection{The term involving $v_\nu$}
We start by considering the term of $e_\nu$ depending on $v_\nu$.
Since
$$
\partial _t v_\nu = T_{\alpha_\nu^{-1/2}}\partial^2_t u_\nu - T_{\partial^2_t(\alpha_\nu^{-1/2})}u_\nu\,,
$$
using also equation \eqref{eq:u_nu} we deduce that
\begin{align*}
\frac{d}{dt}\left\|v_\nu(t,\cdot )\right\|^2_{L^2}
&=2\,\Re\big( v_\nu(t,\cdot ), \partial _t v_\nu (t,\cdot )\big)_{L^2} \\[0.2cm]
&=2\,{\Re}\left( v_\nu, T_{\alpha_\nu^{-1/2}}\partial^2_t u_\nu\right)_{L^2} -
2\,{\Re} \left( v_\nu, T_{\partial^2_t(\alpha_\nu^{-1/2})}u_\nu\right)_{L^2}
\\[0.2cm]
&= -2\,{\Re} \left( v_\nu, T_{\partial^2_t(\alpha_\nu^{-1/2})}u_\nu\right)_{L^2}
+ 2\,{\Re} \left(v_\nu,T_{\alpha_\nu^{-1/2}}\left( \sum_{j,k}\partial_j\left(T_{a_{jk}}\partial_k u_\nu\right)\right)\right)_{L^2} \\[0.2cm]
&\phantom{=}\; + 2\,{\Re} \left(v_\nu, \sum_{j,k}T_{\alpha_\nu^{-1/2}}\partial_j\Big(\left[\Delta_\nu,\,T_{a_{jk}}\right]\partial_k u\Big)\right)_{L^2} +
2\,{\Re} \left(v_\nu, T_{\alpha_\nu^{-1/2}}\left(\wtilde L u\right)_\nu\right)_{L^2}\,.
\end{align*}

Observe that
\begin{align*}
\left|{\Re} \left(v_\nu, T_{\alpha_\nu^{-1/2}}\left(\wtilde L u\right)_\nu\right)_{L^2}\right| &\leq 
C\,\big(e_\nu(t)\big)^{\frac{1}{2}}\, \left\|\left(\wtilde L u\right)_\nu\right\|_{L^2}\,. 
\end{align*}
In addition, whenever $\nu\geq1$ Bernstein's inequalities (see Lemma 2.1 in \cite{BCD}) allow us to bound
\begin{align*}
\left|{\Re} \left( v_\nu, T_{\partial^2_t(\alpha_\nu^{-1/2})}u_\nu\right)_{L^2}\right| &\leq C\,
\left\|v_\nu\right\|_{L^2}\,\big(\varphi_\nu(t)\big)^2\,\left\|u_\nu\right\|_{L^2}\,\leq\,C\,\big(\varphi_\nu(t)\big)^2\, 2^{-\nu} \, e_\nu(t)\,;
\end{align*}
remark that, the previous inequality is trivially verified also in the case $\nu=0$.

In view of the previous bounds, we obtain
\begin{align}
\label{estdtvnu2}
\frac{d}{dt}\left\|v_\nu(t,\cdot )\right\|^2_{L^2}
&\leq2\,{\Re} \left(v_\nu,T_{\alpha_\nu^{-1/2}}\left( \sum_{j,k}\partial_j\left(T_{a_{jk}}\partial_k u_\nu\right)\right)\right)_{L^2} \\[0.2cm]
\nonumber
&\phantom{=}\qquad + 2\,{\Re} \left(v_\nu, \sum_{j,k}T_{\alpha_\nu^{-1/2}}\partial_j\Big(\left[\Delta_\nu,\,T_{a_{jk}}\right]\partial_k u\Big)\right)_{L^2} \\[.2cm]
\nonumber
&\phantom{=}\qquad\qquad + C\,\big(e_\nu(t)\big)^{\frac{1}{2}}\, \left\|\left(\wtilde L u\right)_\nu\right\|_{L^2}\,+\,C\,\big(\varphi_\nu(t)\big)^2\, 2^{-\nu} \, e_\nu(t)\,.
\end{align}

\subsubsection{Estimate of the term involving $w_\nu$}
We now compute the time derivative of the last term appearing in the definition of $e_\nu$, namely the one depending on $w_\nu$. An easy computation
shows that
$$
\frac{d}{dt}\left\|w_\nu(t,\cdot)\right\|_{L^2}^2= 2\,{\Re}\left(T_{\partial_t(\alpha_\nu^{1/2})(\gamma^2+|\xi|^2)^{1/2}} u_\nu + 
T_{\alpha_\nu^{1/2}(\gamma^2+|\xi|^2)^{1/2}} \partial_t u_\nu, w_\nu (t,\cdot)\right)_{L^2}\,.
$$
We claim that
\begin{equation}\label{estdtwnu2}
\frac{d}{dt}\left\|w_\nu(t,\cdot)\right\|^2_{L^2}\,=\, 2\,{\Re}\left(v_\nu , T_{\alpha_\nu^{-1/2}} T_{\alpha_\nu^2(\gamma^2+|\xi|^2)}u_\nu\right)_{L^2}\,+\,Q\,,
\end{equation}
where the remainder term $Q$ satisfies the bound
\begin{equation} \label{est:claimed}
\left|Q\right|\,\leq\,C\,\left(\nu^\ell+1\right)\,e_\nu(t)\,,
\end{equation}
where $\ell\in\{0,1\}$ has been introduced in Section \ref{s:results}.

In order to prove our claim, we follow the computations of \cite{CDSFM1}, \cite{CDSFM2} and freely use the facts proven therein,
about the principal symbols of the remainder operators and their orders. We start by writing
\begin{align}
\label{estdtwnu12} 
&{\Re}\left(T_{\partial_t(\alpha_\nu^{1/2})(\gamma^2+|\xi|^2)^{1/2}} u_\nu, w_\nu\right)_{L^2}\\[0.2cm]
\nonumber
&\qquad\qquad\qquad
={\Re}\left( T_{\alpha_\nu(\gamma^2+|\xi|^2)^{1/2}} T_{-\partial_t(\alpha_\nu^{-1/2})}u_\nu,w_\nu\right)_{L^2} +{\Re}\left(R_1 u_\nu, w_\nu\right)_{L^2}\,,
\end{align}
where the principal symbol of $R_1$ is
$$
\partial_\xi\left(\alpha_\nu (\gamma^2+|\xi|^2)^{1/2}\right)\partial_x\partial_t\left(\alpha_\nu^{-1/2}\right)
$$
and consequently one has
\begin{align*}
\big|{\Re} \left(R_1 u_\nu, w_\nu\right)_{L^2}\big|\,&\leq\, C\, \varphi_\nu(t)\, \left(1+\nu^\ell\right)\, \left\|u_\nu\right\|_{L^2}\,\left\|w_\nu\right\|_{L^2}\,
\leq\, C\, \varphi_\nu(t)\, \left(1+\nu^\ell\right)  \,2^{-\nu}\,\left\|w_\nu\right\|_{L^2}^2 \\
&\leq\, C\,\left(\nu^\ell+1\right)\, e_\nu(t)\,,
\end{align*}
where we have used also estimate \eqref{boundphi2}.
Observe that passing from the first to the second inequality relies on the use of the Bernstein inequality, thus requires $\nu\geq1$; however, for $\nu=0$
one can directly pass from the first inequality to the last one.

Next, we can compute
\begin{align}
\label{estdtwnu22} 
&{\Re} \left(T_{\alpha_\nu^{1/2}(\gamma^2+|\xi|^2)^{1/2}} \partial_t u_\nu, w_\nu\right)_{L^2} \\[.2cm]
\nonumber
&\qquad\qquad\qquad =\,
{\Re} \left( T_{\alpha_\nu(\gamma^2+|\xi|^2)^{1/2}} T_{\alpha_\nu^{-1/2}} \partial_t u_\nu,w_\nu\right)_{L^2}+
{\Re}\left(R_2  \partial_t u_\nu, w_\nu\right)_{L^2}\,,
\end{align}
where, this time, the principal symbol of the remainder $R_2$ is
$$
\partial_\xi\left(\alpha_\nu\left(\gamma^2+|\xi|^2\right)^{1/2}\right)\partial_x\left(\alpha_\nu^{-1/2}\right)\,.
$$
With the help of \eqref{estpartunu2}, we can bound
$$
\big|{\Re}\left(R_2  \partial_t u_\nu, w_\nu\right)_{L^2}\big|\,\leq\, C\, \left(1+\nu^\ell\right)\,\left\|\partial_t u_\nu\right\|_{L^2}\,
\left\|u_\nu\right\|_{L^2}\,\leq\, C\, \left(\nu^\ell+1\right)\, e_\nu(t)\,.
$$

Thus, putting together \eqref{estdtwnu12} and \eqref{estdtwnu22}, we find
$$
\frac{d}{dt}\left\|w_\nu(t,\cdot)\right\|^2_{L^2}\,=\,2\,{\Re}\left(T_{\alpha_\nu(\gamma^2+|\xi|^2)^{1/2}}\, v_\nu, w_\nu\right)_{L^2}\,+\,\wtilde Q\,,
$$
where the remainder $\wtilde Q$ satisfies the claimed bound \eqref{est:claimed}.
In addition, we can write
$$
{\Re}\left(T_{\alpha_\nu(\gamma^2+|\xi|^2)^{1/2}}\, v_\nu, w_\nu\right)_{L^2}=
{\Re}\left( v_\nu,T_{\alpha_\nu(\gamma^2+|\xi|^2)^{1/2}} w_\nu\right)_{L^2} +{\Re}\left(R_3 v_\nu, w_\nu\right)_{L^2}\,,
$$
where the principal symbol of $R_3$ is given by
$$
\partial_\xi \partial_x\left(\alpha_\nu\left(\gamma^2+|\xi|^2\right)^{1/2}\right)\,.
$$
Consequently, keeping in mind the estimates of Subsection \ref{ss:symbolic}, we see that $R_3$ is an operator of order $0$ in the case $\ell=0$,
of order ${\rm Log}$ if $\ell=1$. This means that
\begin{equation} \label{est:R_3}
\big|{\Re}\left(R_3v_\nu, w_\nu\right)_{L^2}\big|\,\leq\, C\,\left(1+\nu^\ell\right)\,\left\|v_\nu\right\|_{L^2}\,\left\|w_\nu\right\|_{L^2}\,\leq\,
C\, \left(\nu^\ell+1\right)\, e_\nu(t)\,.
\end{equation}

Next, we continue to decompose
$$
{\Re}\left( v_\nu,T_{\alpha_\nu(\gamma^2+|\xi|^2)^{1/2}} w_\nu\right)_{L^2}
={\Re} \left( v_\nu,T_{\alpha_\nu^{-1/2}}T_{\alpha_\nu^{3/2}(\gamma^2+|\xi|^2)^{1/2}} w_\nu\right)_{L^2} +{\Re}\left(R_4  v_\nu, w_\nu\right)_{L^2}\,,
$$
where, by using again symbolic calculus, it is plain to see that also the new remainder $R_4$ is an operator of order ${\rm Log}$, in the sense
that is satisfies an estimate like \eqref{est:R_3}.

Finally, we can write
$$
{\Re} \left( v_\nu,T_{\alpha_\nu^{-1/2}}T_{\alpha_\nu^{3/2}(\gamma^2+|\xi|^2)^{1/2}} w_\nu\right)_{L^2}
={\Re}\left( v_\nu,T_{\alpha_\nu^{-1/2}}T_{\alpha_\nu^{2}(\gamma^2+|\xi|^2)} u_\nu\right)_{L^2} +{\Re} \left(R_5 v_\nu, u_\nu\right)_{L^2}\,,
$$
where the principal symbol of $R_5$ is
$$
\partial_\xi \left(\alpha_\nu^{3/2}\left(\gamma^2+|\xi|^2\right)^{1/2}\right)\,\partial_x \left(\alpha_\nu^{-1/2}\left(\gamma^2+|\xi|^2\right)^{1/2}\right)\,.
$$
The previous expression implies that $R_5$ is an operator of order $1$ if $\ell=0$, of order $1+{\rm Log}$ when $\ell=1$; using the Bernstein inequalities
to absorbe the additional power of $\nu$ coming out in the estimate (with the usual modification when $\nu=0$ for dealing
with the low frequency term), we can bound
\begin{align*}
\big|{\Re}\left(R_5v_\nu, u_\nu\right)_{L^2}\big|\,&\leq\, C\, 2^\nu \left(1+\nu^\ell\right)\,\left\|v_\nu\right\|_{L^2}\,\left\|u_\nu\right\|_{L^2}\,\leq\,
 C\,\left(1+\nu^\ell\right)\,\left\|v_\nu\right\|_{L^2}\,\left\|w_\nu\right\|_{L^2} \\
 &\leq\, C\,\left(\nu^\ell+1\right)\, e_\nu(t)\,.
\end{align*}
In the end, the claimed relations \eqref{estdtwnu2}-\eqref{est:claimed} follow.

\subsubsection{Summing up $v_\nu$ and $w_\nu$}

At this point, we would like to put formulas \eqref{estdtunu2}, \eqref{estdtvnu2}, \eqref{estdtwnu2} and \eqref{est:claimed} together.
For this, we first need to handle the sum
$$
{\Re} \left(v_\nu, T_{\alpha_\nu^{-1/2}}\left(\sum_{j,k}  \partial_j\left(T_{a_{jk}}\partial_k u_\nu\right)\right)\right)_{L^2} +
{\Re} \left(v_\nu , T_{\alpha_\nu^{-1/2}} T_{\alpha_\nu^2(\gamma^2+|\xi|^2)}u_\nu\right)_{L^2}\,.
$$

To begin with, we remark that we can write
\begin{equation} \label{eq:ellipt}
\sum_{j,k} \partial_j\left(T_{a_{jk}}\partial_k u_\nu\right)= \sum_{j,k} R_{jk} u_\nu - T_{\sum_{j,k}{a_{jk}\xi_j\xi_k }}u_\nu\,,
\end{equation}
where, recalling the definition \eqref{def:paraprod} of the paraproduct operator, we have set
$$
R_{jk} u_\nu = \partial_j\left( S_{\mu -1}a_{jk}\right)\,S_{\mu +2}\partial_k u_\nu + \sum_{h=\mu-3}^{+\infty} \partial_j\left( S_{h-3} a_{jk}\right)
\Delta_{h}\partial_k u_\nu\,.
$$
Recall that $\mu=\mu(\g)$ is now fixed, owing to our choice of $\g$. Consequently, as done in Paragraph 4.2.4 of \cite{CDSFM1}, we can bound
$$
\left\|R_{jk} u_\nu\right\|_{L^2}\,\leq\, C\, \left(\nu^\ell+1\right)\, e_\nu(t)\,,
$$
for a constant $C>0$ depending on the chosen $\gamma$ and on the Lipschitz norms (if $\ell=0$) or log-Lipschitz norms (when $\ell=1$) of the $a_{jk}$'s.

On the other hand, by definition of $\alpha_\nu$, we have
$$
T_{\alpha_\nu^2(\gamma^2+|\xi|^2)}u_\nu\,=\, T_{\gamma^2 +\sum_{j,k}a_{\nu, jk}\xi_j\xi_k }u_\nu\,.
$$
From this expression, we deduce the equality
$$
T_{\sum_{j,k}a_{\nu, jk}\xi_j\xi_k }u_\nu-T_{\sum_{j,k}{a_{jk}\xi_j\xi_k }}u_\nu= T_{\sum_{j,k}(a_{\nu, jk}-a_{jk})\xi_j\xi_k }u_\nu\,.
$$
Observe that, owing to the properties stated in Subsection \ref{ss:time}, one has
\begin{align*}
a_{\nu, jk}-a_{ jk} = 0 \qquad\qquad &\text{ if }\qquad t\geq 2^{-\nu+1}\,,\\[0.2 cm]
\big|a_{\nu, jk}-a_{ jk}\big|\leq C\qquad\qquad&\text{ if }\qquad t\leq 2^{-\nu+1}\,.
\end{align*}
Thanks to these estimates, we gather (by Bernsteing inequality when $\nu\geq1$, by simple computation when $\nu=0$) the bound
\begin{align*}
\left\|T_{\sum_{j,k}(a_{\nu, jk}-a_{jk})\xi_j\xi_k }u_\nu\right\|_{L^2}\,&\leq\,
C\, 2^{2\nu}\,\chi_{[0,2^{-\nu+1}]}(t)\,\left\|u_\nu\right\|_{L^2} \\[0.2cm]
&\leq \,C\, 2^{\nu}\,\chi_{[0,2^{-\nu+1}]}(t)\,
\big(e_\nu(t)\big)^{1/2}\,.
\end{align*}

Consequently, bearing also \eqref{eq:ellipt} in mind, we infer
\begin{align}\label{estdtvnu+wnu2}
&\left|{\Re} \left(v_\nu, T_{\alpha_\nu^{-1/2}}\left(\sum_{j,k}  \partial_j\left(T_{a_{jk}}\partial_k u_\nu\right)+
T_{\alpha_\nu^2(\gamma^2+|\xi|^2)}u_\nu\right)\right)_{L^2}\right| \\
\nonumber
&\qquad\qquad\qquad\qquad\qquad\qquad\qquad\qquad
\leq\, C\,\left(1+\nu^\ell+2^{\nu}\chi_{[0,2^{-\nu+1}]}(t)\right)\,e_\nu(t)\,,
\end{align}
where the constant $C>0$ depends on $\gamma$, $\Lambda_0$ and on the Lipschitz (for $\ell=0$) or log-Lipschitz (for $\ell=1$) norms of the $a_{jk}$'s.

Summing up, putting estimates \eqref {estdtunu2}, \eqref {estdtvnu2}, \eqref {estdtwnu2}, \eqref{est:claimed} and \eqref{estdtvnu+wnu2} together, we deduce
\begin{align} \label{estdtenu2}
\frac{d}{dt}\,e_\nu(t)\,&\leq\,{K}_1\,\left(1+\nu^\ell+\big(\varphi_\nu(t)\big)^2\, 2^{-\nu}+2^{\nu}\chi_{[0,2^{-\nu+1}]}(t)\right)\,e_\nu(t) \\
\nonumber
&\phantom{=}\quad + 2\,{\Re} \left(v_\nu, \sum_{j,k}T_{\alpha_\nu^{-1/2}}\partial_j\Big(\left[\Delta_\nu,\,T_{a_{jk}}\right]\partial_k u\Big)\right)_{L^2} 
\,+ \wtilde{K}_1\,\big(e_\nu(t)\big)^{\frac{1}{2}}\, \left\|\left(\wtilde L u\right)_\nu\right\|_{L^2}\,,
\end{align}
for two suitable positive constants $K_1$ and $\wtilde K_1$, only depending on $\g$, $\lambda_0$, $\Lambda_0$ and on the (Lipschitz of log-Lipschitz) functional norms
of the coefficients $a_{jk}$ of $L$.

\subsection{Estimate of the commutator term}
In order to find an estimate on the time derivative of the localised energy $e_\nu$, it remains us to exhibit a control for the commutator term
appearing in \eqref{estdtenu2} above.
However, it turns out that it is better to handle the whole sum over $\nu$ of all the commutator terms, instead of considering them one by one at $\nu$ fixed.

Here the computation is exactly the same as in \cite{CDSF} if $\ell=0$, as in \cite{CDSFM1}
in the case $\ell=1$. Indeed, no derivatives with respect to time are involved and only the
(Lipschitz or log-Lipschitz) regularity with respect to $x$ of the coefficients is needed.

\medbreak
Let us consider the case $\ell=1$ for a while, and take a $\theta\in \, ]0,1[\,$. Then, for $\beta>0$ and $T_0\in \, ]0,T[\,$ to be chosen later,
we obtain, for all $t\in [0,\, T_0]$, the inequality
\begin{align}\label{estcommLL}
&\sum_{\nu=0}^{+\infty} e^{-2\beta(\nu+1)t}\, 2^{-2\nu\theta}
\left|2\,{\Re} \left(v_\nu, \sum_{j,k}T_{\alpha_\nu^{-1/2}}\partial_j\Big(\left[\Delta_\nu,\,T_{a_{jk}}\right]\partial_k u\Big)\right)_{L^2}\right|\\[0.2 cm]
\nonumber
&\qquad\qquad\qquad\qquad\qquad\qquad\qquad\qquad\qquad
\leq\,K_2\,\sum_{\nu=0}^{+\infty}(\nu+1)\,e^{-2\beta(\nu+1)t}\, 2^{-2\nu\theta} e_\nu(t)\,,
\end{align}
where the constant $K_2>0$ depends on $\gamma$, $\theta$, the product $\beta T_0$  and the log-Lipschitz norm of the $a_{jk}$'s.
We refer to Subsection 4.3 of \cite{CDSFM1} for details.

\medbreak
Consider now the case $\ell=0$. Take any $\theta\in[0,1[\,$. Then, according to estimate (37) in \cite{CDSF}, for any $t\in[0,T]$ one has instead
\begin{align}\label{estcommLip}
&\sum_{\nu=0}^{+\infty} 2^{-2\nu\theta}
\left|2\,{\Re} \left(v_\nu, \sum_{j,k}T_{\alpha_\nu^{-1/2}}\partial_j\Big(\left[\Delta_\nu,\,T_{a_{jk}}\right]\partial_k u\Big)\right)_{L^2}\right|\,
\leq\,K_2\,\sum_{\nu=0}^{+\infty} 2^{-2\nu\theta} e_\nu(t)\,,
\end{align}
for a new positive constant, that we keep calling $K_2$ for convenience. Notice that, as above, $K_2$ only depends on
$\g$, $\theta$ and on the Lipschitz norms of the coefficients $a_{jk}$.

\subsection{The total energy: end of the estimates}

It is now time to 
define the total energy $\mc E$ related to $u$ and close the estimate for such $\mc E$.

First of all, let us define the function $f_\nu(t)$, for any $t\in[0,T_0]$ (where we agree that $T_0\equiv T$ when $\ell=0$), as
\[
 f_\nu(t)\,:=\,\int^t_0\left(\big(\varphi_\nu(\t)\big)^2\, 2^{-\nu}+2^{\nu}\chi_{[0,2^{-\nu+1}]}(\t)\right)\,d\t\,.
\]
Observe that, for all $t\in [0,\,T_0]$, one has
\begin{equation}\label{estint2}
0\,\leq\, f_\nu (t)\, \leq\, C\,,
\end{equation}
where $C>0$ does not depend on $\nu\in\N$. Consequently, when $\ell=1$, from \eqref{estcommLL} we infer 
\begin{align}\label{estcomm3}
&\sum_{\nu=0}^{+\infty} e^{-K_1\,f_\nu(t)}\,e^{-2\beta(\nu+1)t}\, 2^{-2\nu\theta}\,
\left|2\,{\Re} \left(v_\nu, \sum_{j,k}T_{\alpha_\nu^{-1/2}}\partial_j\Big(\left[\Delta_\nu,\,T_{a_{jk}}\right]\partial_k u\Big)\right)_{L^2}\right|\\[0.2 cm]
\nonumber
&\qquad\qquad\qquad\qquad\qquad\qquad\qquad\qquad
\leq\, K_3\, \sum_{\nu=0}^{+\infty} (\nu+1)\,e^{-K_1f_\nu(t)}\, e^{-2\beta(\nu+1)t}\, 2^{-2\nu\theta}\, e_\nu(t)\,,
\end{align}
where $K_3$ depends on the same quantities as $K_2$ depends on, namely on $\gamma$, $\theta$, the product $\beta T_0$  and the log-Lipschitz norm of the $a_{jk}$'s.

When $\ell=0$, instead, we use \eqref{estcommLip} and get
\begin{align}\label{estcomm-Lip_bis}
&\sum_{\nu=0}^{+\infty} e^{-K_1\,f_\nu(t)}\, 2^{-2\nu\theta}\,
\left|2\,{\Re} \left(v_\nu, \sum_{j,k}T_{\alpha_\nu^{-1/2}}\partial_j\Big(\left[\Delta_\nu,\,T_{a_{jk}}\right]\partial_k u\Big)\right)_{L^2}\right|\,
\leq\,K_3\,\sum_{\nu=0}^{+\infty} 2^{-2\nu\theta} e_\nu(t)\,,
\end{align}
where in fact $K_3$ can be taken exactly equal to the constant $K_2$ appearing in \eqref{estcommLip}.

We are now in the position of defining the total energy $\mc E$ and closing the estimates. 
From now on, we will keep our argument general, without specifing each time whether $\ell=0$ or $\ell=1$. However, we will tacitly agree on the following conditions:
\begin{itemize}
 \item if $\ell=1$, we take $\theta\in\,]0,1[\,$, $\beta>0$ and $0<T_0\leq T$, where both $\beta$ and $T_0$ have to be fixed later;
 \item if instead $\ell=0$, we take $\theta\in[0,1[\,$, $\beta=0$ and $T_0=T$.
\end{itemize}

This having been pointed out, let us define
$$
{\cal E}(t)\,:=\, \sum_{\nu=0}^\infty e^{-K_1 f_\nu(t)}  \,  e^{-2\beta(\nu+1)t}\, 2^{-2\nu\theta} e_\nu(t)\,.
$$
Then, combining \eqref {estdtenu2} with  \eqref {estcomm3} when $\ell=1$, or with \eqref{estcomm-Lip_bis} for $\ell=0$, we deduce that 
\begin{align*}
{\cal E}'(t)\, &\leq\,\big(K_1+ K_3-2\beta\big)\sum_{\nu=0}^{+\infty} (\nu+1)\,e^{-K_1  f_\nu(t)}\,  e^{-2\beta(\nu+1)t}\, 2^{-2\nu\theta} e_\nu(t) \\
&\qquad\qquad\qquad\qquad
+\,K_4\,  \sum_{\nu=0}^{+\infty} 
e^{-K_1  f_\nu(t)}\,  e^{-2\beta(\nu+1)t}\, 2^{-2\nu\theta}\, \big(e_\nu(t)\big)^{\frac{1}{2}}\,\left\|\left(\wtilde L u\right)_\nu\right\|_{L^2}\,,
\end{align*}
for another ``universal'' constant $K_4>0$.

At this point, the argument for bounding the term $\wtilde Lu$ is the same as \cite{CDSFM1} (for $\ell=1$) and \cite{CDSF} (for $\ell=0$),
the only difference being the presence of the term $e^{-K_1  f_\nu(t)}$ which however, in view of \eqref{estint2}, plays the role of a constant.
Thus, for $\ell=1$ we obtain that, for any $\theta\in \, ]0,1[\,$, there exist $\beta^*\,=\,\beta\,(\log2)^{-1}$, $T_0\in\,]0,T]$ and $C>0$ such that
$$
{\cal E}'(t)\,\leq\, C\,\left({\cal E}(t)\, +\,\big({\cal E}(t)\big)^{1/2}\, \left\|Lu(t, \cdot)\right\|_{H^{-\theta-\beta^*t}}\right)
$$
for all $u\in C^2\big([0,T_0]; H^\infty(\R^n)\big)$. Once again, if $\ell=0$, then $\theta$ can take the value $0$ and the previous estimate holds with $\beta^*=0$
and $T_0=T$.

In the end, the claimed estimate follows from the fact that
\begin{align*}
{\cal E}(0)\,&\leq\, C_\theta\, \left(\left\| u(0,\cdot)\right\|^2_{H^{1-\theta}}\,+\,\left\|\partial_tu(0,\cdot)\right\|^2_{H^{-\theta}}\right)\,,\\[0.2 cm]
{\cal E}(t)\,&\geq\, C'_\theta\, \left(\left\| u(t,\cdot)\right\|^2_{H^{1-\theta-\beta^*t}}\,+\,\left\|\partial_tu(t,\cdot)\right\|^2_{H^{-\theta-\beta^*t}}\right)\,.
\end{align*}
Observe that, in fact, the constant $C'_\theta$ also depends on $K_1$, hence on the quantities $\g$, $\lambda_0$, $\Lambda_0$ and on the functional
norms of the coefficients $a_{jk}$, namely on the constants $C_0\ldots C_3$ appearing in assumptions (H1) and (H2) of Section \ref{s:results}.

This concludes the proof of Theorems \ref{th:Yam} and \ref{theorem 3}.


\end{document}